\documentclass[b5paper,10pt,twoside,onecolumn,draft]{article}
\begin{document}


\title{Fast Computation of the Arnold Complexity\protect\\
of Length $2^{n}$ Binary Words}

\author{Yuri V. Merekin}

\maketitle

 {\small Sobolev Institute of Mathematics, SB RAS,
Novosibirsk 630090,
Russia. \\
E-mail: merekin@math.nsc.ru}



{\bf Abstract.}{\small\ For fast computation of the Arnold
complexity of length $2^{n}$ binary words we obtain an upper bound
for the Shannon function $Sh(n)$.


\section{Introduction}

Analyzing word complexity usually involves studying the fragments of
a~word or the process of its construction (see~[2] for instance).
Arnold introduced~[1] a~new concept of complexity of a~word. The
measure of this complexity is determined by the ``stability'' of
a~word under the iterated action of a~certain operator.


Consider an arbitrary nonperiodic binary word $w=x_{1}x_{2}\ldots
x_{2^{n}}$, with $w\neq v^{k}$ and $k\geq 2$, of length $|w|=2^{n}$
for $n\geq 1$. Denote by $(w)$ the infinite periodic word $(w)=ww
\ldots$. Henceforth,
by a ``word~%
$(w)$''
we understand
``an infinite periodic word~%
$(w)$''.
Consider the scheme
(word chain)
\begin{equation}                                                  
(w)=(w_{1}), (w_{2}), \ \ldots \ , (w_{s})=(v),
\end{equation}
in which the first word is arbitrary,
and every word
$(w_{i})=y_{1}y_{2} \ldots$
generates the next word
$(w_{i+1})=z_{1}z_{2} \ldots$,
for
$1\leq i\leq s-1$,
using the operator
\begin{equation}                                                  
F(\cdot , h_{i}):(w_{i})\mapsto (w_{i+1}): \ z_{j}=y_{j} \oplus y_{j+h_{i}},
\end{equation}
where
$j\geq 1$,
$1\leq h_{i}=2^{n_{i}}$,
$0\leq n_{i}\leq n$,
and~%
$\oplus$
stands for modulo~2 addition;
thus,
$F(w_{i},h_{i})=(w_{i+1})$.
The number
$h_{i}$
is called the \textit{rank} of the operator in~(2),
and the number~%
$s$,
the \textit{length} of the scheme~(1).

Denote  by
$S(w,v,h,s)$
the type of schemes with the first word
$(w)$,
the last word
$(v)$,
the maximal rank~%
$h$
of operators involved,
and the scheme length~%
$s$.
For every
$(w)$
there exists a minimal~%
$s$
such that
all words in the scheme
$S(w,0,1,s)$
are distinct,
and
$F(v,1)=(0)$.
A scheme of this type is called a~\textit{complexity scheme}.
Every word
$(w)$
has a unique complexity scheme.
The number
$s-1$
is called the \textit{complexity } of the binary word
$(w)$
and is denoted by
$A(w)$.
%
%
Arnold introduced~[1]
the concept of complexity of a binary word
in a more general form,
which coincides with our definition of complexity
when the word length equals~%
$2^{n}$.
%
%
The complexity of a periodic word
$(w)$
is equal to the complexity of the finite word~%
$w$.

In an arbitrary scheme
$S(w,v,1,s)$
select a word chain
\begin{equation}                                                 
(w)=(w_{1}), (w_{1+l_{1}}), \ \ldots \ , (w_{1+l_{1}+ \ldots
+l_{t}})=(v),
\end{equation}
where
$l_{i}\geq 1$
for
$1\leq i\leq t\leq s-1$.
If in~(3) each word
$(w_{1+l_{1}+ \ldots +l_{i}})$
for
$1\leq i\leq s-1$,
coincides with
$F(w_{1+l_{1}+ \ldots +l_{i-1}},h_{i})$,
$h_{i}=l_{i}$,
then the word chain in~(3) is a scheme of type
$S(w,v,h,s_{t})$
with
$s_{t}=1+t$,
which is called \textit{equivalent} to 
$S(w,v,1,s)$.

In~[3] we proved

{\bf Theorem 1.1} {\it Every scheme $S(w,v,1,2^{n}+1)$ with $n \geq
0$ is equivalent to the elementary scheme $S(w,v,2^{n},2)$.}

In a word
$w=x_{1}x_{2}\ldots x_{2^{n}}$,
$n\geq1$,
select
$2^{n-m}$
positions,
with
$0\leq m\leq n-1$,
such that
the distances between two neighboring selected positions
are the same and equal to~%
$2^{m}$.
Using the selected positions,
form the word
$u=x_{i}x_{i+2^{m}} \ldots x_{i+2^{n}-2^{m}}$
of length
$2^{n-m}$.
Denote the infinite word
$(u)=uu\ldots$
by
$(w)^{x_{i}}_{2^{n-m}}$
and call it a~\textit{thinned-out word}.
The number
$2^{n-m}$
is called the \textit{step} of the thinned-out word
$(w)^{x_{i}}_{2^{n-m}}$.
Observe that
every thinned-out word
is a linearly ordered set of indices of positions of~%
$w$.
The length of the period of
$(w)^{x_{i}}_{2^{n-m}}$
can be less than
$2^{n-m}$.
For
$m=0$
we have
$(w)=(w)^{x_{1}}_{2^{n}}$.

Given a word
$(w)$,
for a fixed value of~%
$m$
there exist
$2^{m}$
different thinned-out words
\begin{equation}                                                 
(w)^{x_{1}}_{2^{n-m}} \ , (w)^{x_{2}}_{2^{n-m}} \ , \ \ldots  \ ,
(w)^{x_{2^{m}}}_{2^{n-m}}, \ \ n\geq1, \ \ 0\leq m\leq n-1.
\end{equation}
For instance, for $n=3$ we have $
(w)^{x_{1}}_{2^{3}}=(x_{1}x_{2}x_{3}x_{4}x_{5}x_{6}x_{7}x_{8}), m=0;
\\
(w)^{x_{1}}_{2^{3-1}}=(x_{1}x_{3}x_{5}x_{7}),\
(w)^{x_{2}}_{2^{3-1}}=(x_{2}x_{4}x_{6}x_{8}), \ m=1;
\\
(w)^{x_{1}}_{2^{3-2}}=(x_{1}x_{5}), \
(w)^{x_{2}}_{2^{3-2}}=(x_{2}x_{6}), \
(w)^{x_{3}}_{2^{3-2}}=(x_{3}x_{7}),\\
\hspace*{7.5cm}(w)^{x_{4}}_{2^{3-2}}=(x_{4}x_{8}),  m=2. $

Define the operation of taking the union of thinned-out words,
denoted by the symbol~%
$\ast$.
The definition of a thinned-out word implies that
each of the positions
$x_{1}, x_{2}, \ldots , x_{2^{n}}$
appears in the thinned-out words~(4) exactly once
since it is the union of arithmetic progressions
with differences equal to powers of~2.
Thus,
we can express
$(w)$
as
$$
(w)= (w)^{x_{1}}_{2^{n-m}} \ast (w)^{x_{2}}_{2^{n-m}} \ast \dots
\ast (w)^{x_{2^m}}_{2^{n-m}} \ \ n\geq1, \ \ 0\leq m\leq n-1.
$$

We group the words into two sorts:
even words and odd words.
The word
$(w) = ww\ldots$, where
$w=x_{1}x_{2} \ldots x_{2^{n}}$
with
$x_{i}\in \{0,1\}$
and
$1\leq i\leq 2^{n}$
for
$n\geq 0$,
is called \textit{even} whenever
$x_{1}\oplus x_{2} \oplus \ldots \oplus x_{2^{n}} = 0$,
and \textit{odd} whenever
$x_{1}\oplus x_{2} \oplus \ldots \oplus x_{2^{n}} = 1$.
%
%
For calculating the parity of the thinned-out words
$(w)^{x_{i}}_{2^{n-m}}$, where $|w|=2^{n}$ with $n\geq1$ and $0\leq
m\leq n-1$, of a word $(x_{1}x_{2}\ldots x_{2^{n}})$, we gave a
simple algorithm~[4],
which uses 
modulo~2 addition $2^{n}-1$ times, and proved

{\bf Theorem1.2} {\it For every binary word $(w)$ the length of
whose period is equal to $2^{n}$, $n\geq1$, all thinned-out words
$(w)^{x_{i}}_{2^{n-m}}$ for $0 \leq m \leq n-1$ and $1\leq i \leq
2^{m}$ are odd if and only if the complexity of $(w)$ is equal to
$A(w) = 2^{n}-2^{m}+1.$}

Express the complexity
$A(w)$
of an arbitrary word
$(w)$
with
$|w|=2^{n}$
for
$n\geq 1$
as 
\begin{equation}
A(w)=a_{n-1}2^{(n-1)}+a_{n-2}2^{(n-1)-1}+ \ldots \ + a_{0}2^{(n-1)-(n-1)},      
\end{equation}
or the binary number
$a_{n-1} a_{n-2} \ldots a_{0}$,
where
$a_{i}\in \{0,1\}$
for
$0 \leq i \leq n-1$.

Express the complexities $A(w)$, which according to Theorem~1.2 we
calculate by finding the parities of thinned-out words, as the
binary numbers
\begin{eqnarray}
2^{n}-2^{1}+1&=&a_{n-1}a_{n-2}\ldots a_{0}=11\ldots 111,\nonumber
\\
2^{n}-2^{2}+1&=&a_{n-1}a_{n-2}\ldots a_{0}=11\ldots 101, \ \ldots  \ , \\        
2^{n}-2^{n-1}+1&=&a_{n-1}a_{n-2}\ldots a_{0}=10\ldots 001, \nonumber
\\
2^{n}-2^{n-1}+0&=&a_{n-1}a_{n-2}\ldots a_{0}=10\ldots 000, \ \ \
n\geq 1, \nonumber
\end{eqnarray}
where all numbers are odd
with the exception of
$2^{n}-2^{n-1}$.

Refer to a word $(v)$
as \textit{final} 
if 
$A(v)$
equals one of the values in~(6). 
Every complexity scheme
$S(w,0,1,s)$
with
$s=2^{n}$
for
$n\geq 1$
contains
$n+(n-1)+\ldots +1$
final words.

Using
$t\geq 1$
operators~(2) of ranks
$h_{1}$,
$h_{2}, \ldots$,
$h_{t}$,
transform the complexity scheme
$S(w,0,1,s)$
into a scheme
$$
(w)=(w_{1}), (w_{1+h_{1}}), \ldots , (w_{1+h_{1}+ \ldots
+h_{t}})=(v), (w_{1+h_{1}+ \ldots +h_{t}+1}),  \ldots , (0)
$$
with the final word~%
$(v)$.
Then
\begin{equation}
A(w)=h_{1}+h_{2}+ \ldots \ +h_{t}+A(v),                                        
\end{equation}
where
$A(v)$
is one of the numbers~(6).
It is obvious that
in order to transform 
$(w)$
into 
$(v)$, every permutation of the ranks $h_{1}$, $h_{2}, \ldots$,
$h_{t}$ of the operators $F(\cdot , h_{i})$ for $1\leq i\leq t$ is
admissible.

\section{Shannon Function}

Refer to the minimal number of operators $F(\cdot , h_{i})$
required to transform 
$(w)$
into one of the final words
$(v)$
as the \textit{complexity of transformation} of 
$(w)$
into 
$(v)$,
 and denote it by
$A(w,v)=\min\limits_{w\rightarrow v} A(w)$.

Our goal is to find the Shannon function
$\max\limits_{w}\min\limits_{w\rightarrow v} A(w)$, which we denote
by~$Sh(n)$.

Consider an example. Take a complexity scheme $(w_{16}), (w_{15}),
\ldots , (w_{2}), (w_{1})$, $(w_{0})$. For the words $(w)$ with
$|w|=2^{n}$ for $0\leq n\leq 4$ five values of complexity exist, for
which we have expressions as in~(7). In each of these cases an
operator $F(w_{i}, 2^{r}):(w_{i})\mapsto (w_{j})$ is used only once.
Table~1 presents the results of calculating 
$A(w)$ for all words $(w)$ with $|w|=2^{n}$ for $1\leq n\leq 4$.

${}$\hfill {Table 1}\hspace*{10pt}
\begin{center}\small
\begin{tabular}{|c|l|l|l|}
\hline
$(w_{i})$&$A(w_{i})$&$F(w_{i}, 2^{r}):(w_{i})\mapsto
(w_{j})$&$A(w_{i})=2^{r}+A(w_{j})$
\\ \hline
\hline
$(w_{16})$
&
$2^{4}-2^{0}+1$
&                                             &  \\
\hline
$(w_{15})$
&
$2^{4}-2^{1}+1$
&                                             &  \\
\hline
$(w_{14})$
&
$
$ &
$F(w_{14}, 2^{0}):(w_{14})\mapsto (w_{13})$
&
$A(w_{14})=2^{0}+(2^{4}-2^{2}+1)$
\\
\hline
$(w_{13})$
&
$2^{4}-2^{2}+1$
&                                             &  \\
\hline $(w_{12})$ & $ $ & $F(w_{12},
2^{2}):(w_{12})\mapsto (w_{8})$
& $A(w_{12})=2^{2}+(2^{3}-2^{0}+1)$
\\
\hline
$(w_{11})$
&
$
$ &
$F(w_{11}, 2^{1}):(w_{11})\mapsto (w_{9})$
&
$A(w_{11})=2^{1}+(2^{4}-2^{3}+1)$
\\
\hline
$(w_{10})$
&
$
$ &
$F(w_{10}, 2^{0}):(w_{10})\mapsto (w_{9})$
&
$A(w_{10})=2^{0}+(2^{4}-2^{3}+1)$
\\
\hline
$(w_{9} )$
&
$2^{4}-2^{3}+1$
&                                             &  \\
\hline
$(w_{8} )$
&
$2^{3}-2^{0}+1$
&                                             &  \\
\hline
$(w_{7} )$
&
$2^{3}-2^{1}+1$
&                                             &  \\
\hline
$(w_{6} )$
&
$
$ &
$F(w_{6}, 2^{0}):(w_{6})\mapsto (w_{5})$
&
$A(w_{6})=2^{0}+(2^{3}-2^{2}+1)$
\\
\hline
$(w_{5} )$
&
$2^{3}-2^{2}+1$
&                                             &  \\
\hline
$(w_{4} )$
&
$2^{2}-2^{0}+1$
&                                             &  \\
\hline
$(w_{3} )$
&
$2^{2}-2^{1}+1$
&                                             &  \\
\hline
$(w_{2} )$
&
$2^{1}-2^{0}+1$
&                                             &  \\
\hline
$(w_{1} )$
&
$2^{0}
$
&                                             &  \\
\hline
\end{tabular}
\end{center}

\medskip In the next theorem we consider the general case for~$n\geq
5$.

 {\bf Theorem 2.1} {\it Given a word $(w)$ with $|w|=2^{n}$ for $n\geq 5$,
we have
$$
Sh(n) \leq \left\{
\begin{array}{ll}
\lfloor n-2\sqrt{n}+1\rfloor & \mbox{when the binary number $A(w)$
is odd},
\\
\lfloor n-2\sqrt{n}+2\rfloor & \mbox{otherwise}.
\end{array}
\right.
$$}

{\it Proof.} \textit{Case 1.} Assume that the value of the
complexity $A(w)$ is odd.

Fix 
$A(w)$ and estimate the minimal number of operators
transforming 
$(w)$
into a final word~%
$(v)$. It is obvious that
 in this case the ranks of all operators are distinct.

In order to estimate $A(w,v)$, consider the result of the action of
the operator of~(2) on the coefficients $a_{n-1}, a_{n-2}, \ldots,
a_{0}$ of~(5). if
$$
 F(u_{1},h=2^{(n-1)-i})=(u_{2}), \ \ 0\leq i\leq n-1,
$$
then $A(u_{1})-A(u_{2})=2^{(n-1)-i}$. Moreover, two variants are
possible for changing the values of $a_{n-1}, a_{n-2}, \ldots,
a_{0}$:
\begin{equation}
(a_{(n-1)-i}=1) \ \mapsto \ (a_{(n-1)-i}=0);
\end{equation}
\begin{equation}
\begin{array}{rcl}
(a_{(n-1)-i+j}=1)     & \mapsto & (a_{(n-1)-i+j}=0),  \\
(a_{(n-1)-i+(j-1)}=0) & \mapsto & (a_{(n-1)-i+(j-1)}=1), \ \ldots \ ,  \\
(a_{(n-1)-i}=0)       & \mapsto & (a_{(n-1)-i}=1),                                 
\end{array}
\end{equation}
where $1\leq i\leq n-1, \ 1\leq j\leq i$.

In case~(8) the rank of 
$h=2^{(n-1)-i}$
coincides with one of the terms in the sum~(5).
The action of the operator removes the term
$2^{(n-1)-i}$
from~(5).
For instance,
the operator of rank
$h=2^{1}$
transforms
$A(w)=2^{4}+2^{3}+2^{1}+2^{0}$
into 
$2^{4}+2^{3}+2^{0}$.

In case~(9),
when the rank of 
$h=2^{(n-1)-i}$
is distinct from all terms of~(5),
we remove the term
$2^{(n-1)-i+j}$
with minimal~%
$j$.
Simultaneously,
we add to~(5) the terms
\begin{equation}
2^{(n-1)-i+(j-1)}, \ 2^{(n-1)-i+(j-2)}, \ldots , 2^{(n-1)-i}.                    
\end{equation}
For instance,
the operator of rank
$h=2^{1}$
transforms
$A(w)=2^{4}+2^{2}+2^{0}$
into
$2^{4}+2^{1}+2^{0}$.

Consider the case when
we can apply~(9)
in order to calculate
$A(w,v)$.

Suppose that
(5) includes a run 
$a_{i}=a_{i-1}=\ldots = a_{i-l+1}=1$
of neighboring unit coefficients of maximal length,
where
$i\leq n-1$
and
$i+l-1\geq 1$,
which we denote by
$s(i,l)$.
Several runs of maximal length may exist;
for instance,
$A(w)=2^{5}+2^{4}+2^{2}+2^{1}+2^{0}$
includes two such runs:
$s(5,2)$
and
$s(2,2)$.

For a fixed value
$A(w)$
of complexity
choose a~run
$s(i,l)$
 arbitrarily.
If
$A(w)$
is distinct from~(6)
then the sum in~(5),
in addition to the~%
$l$
terms
$2^{i}, 2^{i-1}, \ldots , 2^{i-l+1}$
and the term
$2^{0}$,
also involves~%
$t$
distinct terms
with
$1\leq t\leq n-l-2$.
Once we remove these~%
$t$
terms,
the remaining sum would coincide with one of the sums in~(6).

To remove~$t$
distinct terms from~(5)
using~(8)
we need~%
$t$
operators
$F(\cdot , h_{i})$,
of distinct ranks
$h_{1}$,
$h_{2},\ldots$,
$h_{t}$.
For instance,
in
$A(w)=2^{5}+2^{4}+2^{2}+2^{1}+2^{0}$
choose a run of neighboring unit coefficients
of  maximal length
$s(5,2)$
and remove the terms
$2^{2}$
and
$2^{1}$.
This yields the sum
$2^{5}+2^{4}+2^{0}$,
which coincides with one of the sums in~(6).

The transformation process
$A(w)\mapsto A(v)$
involves a~unique case when
the replacement of the variant~(8) by the variant~(9),
in which the number of terms increases,
fails to increase the number
$F(\cdot , h_{i})$
of operators in the transformation
$A(w)\mapsto A(v)$.
Moreover,
the form of the final word changes:
it additionally includes all terms of~(10).
This happens when in 
$A(w)=2^{j}+2^{i}+\ldots +2^{0}$
with
$j\geq i+2$
we choose a~run
$s(i,l)$
and apply the operator
$F(w,h=2^{i+1})$.
Then we remove from
$A(w)$
the term
$2^{j}$
and transform the run
$s(i,l)$
into the run
$s(j-1,l+(j-i-1))$.
For instance,
for
$A(w)=2^{5}+2^{3}+2^{2}+2^{0}$
choose the run
$s(3,2)$.
Then the operator
$F(w, h=2^{4})$
transforms
$A(w)=2^{5}+2^{3}+2^{2}+2^{0}$
into 
$2^{4}+2^{3}+2^{2}+2^{0}$,
while the run
$s(3,2)$
goes into 
$s(4,3)$.

Consequently,
for removing~%
$t$
distinct terms from~
(5)
the application of~(9) is not necessary,
and for finding
$A(w,v)$
we may use only the operators resulting in~(8).

{\bf Remark} {\it Every odd binary number $A(w)$ includes the term
$2^{0}$, which we do not remove while constructing $A(w,v)$.
Consequently, the operator of rank $h=1$ is not used while obtaining
$A(w,v)$.}

Denote by
$\nu (N)$
the number of~1's in the binary expression
for a nonnegative integer~%
$N$.
The arguments above imply that
$$
A(w,v)=\min\limits_{w\rightarrow v} A(w)=\nu (A(w))-l-1,
$$
where~%
$l$
is the length of the maximal run
$s(i,l)$.

Let us find
$Sh(n)$
for nonfinal words
$(w)$
with
$|w|=2^{n}$
for
$n\geq 1$.

Construct a continuous function
which,
copying 
the process of removal of the maximal number of~1's
from a binary number
$A(w)$,
determines an upper bound for
$Sh(n)$.

Divide a line of integer length
$n\geq 4$
into
$x$
segments,
with
$2\leq x\leq n/\! \ 2$,
of the same length~%
$n/\! x$.
Keeping one of the segments intact, 
remove the beginning of all other segments
to leave only a finite part of unit length.
Then the total length of the removed segments
is estimated by the convex function
$$
f(x)=(x-1)(n/\! x-1),
$$
which has one extremum.
Find the derivative
$f'(x)$
and set it equal to zero:
$$
f'(x)=n/\! x^{2}-1=0.
$$
This yields 
$x=\sqrt {n}$
and the maximal value attained by the function
$f(x)$,
equal to
$$
f(x=\sqrt {n})=n-2\sqrt {n}+1.
$$

For odd binary numbers $A(w)$  Theorem is proved.

\textit{Case 2.} Assume that the value of the complexity $A(w)$ is
even. Estimate the minimal number of operators~(2) required for
calculating $A(w,v)$.

Suppose that
the length~%
$n$
binary number
$A(w)$
includes
$\nu (A(w))$
digits~1,
where
$2\leq \nu (A(w))\leq n-1$,
and
\begin{equation}
a_{(n-1)-j}=1, \ a_{(n-1)-j-1}=0, \ \ldots \ , \ a_{0}=0, \ \ 1\leq j\leq n-2.  
\end{equation}

Two variants for calculating
$A(w,v)$
are possible.

\textit{Subcase 2.1.} From the binary number $A(w)=a_{n-1} a_{n-2}
\ldots a_{0}$, which contains $\nu (A(w))$ digits~1,
 remove
$\nu (A(w))-1$
digits~1
using~(8).
This yields a~binary number
$A(v)$
with a~unique digit~1.
The number of operators
$A(w)\mapsto A(v)$
equals
\begin{equation}
 \nu (A(w))-1.                                                                  
\end{equation}

\textit{Subcase 2.2.} Apply the operator $F(w=w_{1},h=1) \colon
(w_{1})\mapsto (w_{2})$ once. As a result, the even number
$A(w_{1})$ goes into the odd number $A(w_{2})$, and
$A(w_{2})=A(w_{1})-1$. Carry out further calculations according to
the algorithm of case~1, in which by Remark~2.2 the operator~(2) of
rank~1 is not used.

Upon the application of 
$F(w=w_{1} , h=1)$
to the number
$A(w)$
all binary digits in~(11) switch their values
in accordance with~(9).
Therefore,
the number
$A(w_{2})$
includes
$$
\nu(A(w_{2}))=\nu(A(w))-1+j
$$
digits~1.
Removing from
$A(w_{2})$
all digits~1
except for~%
$l$
of those in
$s(i, l)$
and
$a_{0}=1$,
we obtain
$A(v)$
with
$\nu (A(v))=l+1$
digits~1.
The number of operators transforming
$A(w)$
into
$A(v)$
equals
\begin{equation}
\nu (A(w_{2}))-\nu (A(v))+1=\nu(A(w))+j-l-1.                                    
\end{equation}

In order to estimate the complexity
$$
A(w,v)=\min\limits_{w\rightarrow v} A(w)
$$
we choose the variant
with the minimal number of operators.
A comparison of (12) and (13) shows that
 this number occurs in subcase~2.1 for
$j\geq l$
and in subcase~2.2 for
$j\leq l$.

For instance,
for
$A(w)=101110100$
we choose subcase~2.2:
$$
\begin{array}{rcll}
A(w)                   &=& 101110100, & \nu (A(w))=5, \ \ \ \ \ \ \ \ \ \ \ \ j=2;  \\
A(w_{2})               &=& 101110011, & \nu (A(w))-1+j=6, \ \ l=3;  \\
A(v)                   &=& 001110001, & \nu (A(v))=l+1=4;  \\
A(w)\mapsto A(v)&\Longrightarrow & 1\phantom{00000}11\phantom{0}, &
\nu(A(w))+j-l-1=3.
\end{array}
$$

Observe that
if the operator
$F(w=w_{1} , h=1):(w_{1})\mapsto (w_{2})$
in subcase~2.2 generates an odd word
$(w_{2})$,
for which we have already established the estimate
$\lfloor n-2\sqrt{n}+1\rfloor$,
then for even
$A(w_{1})$
we have the estimate
$\lfloor n-2\sqrt{n}+2\rfloor$
since
$A(w_{1})=A(w_{2})+1$.

Supported in part by RFBR grant 11-01-00997.

\end{document}